\documentclass{amsart}
\usepackage{graphicx}
\numberwithin{equation}{section}

\begin{document}
\title[Beyond the PI Controllers in First-Order Time-Delay Systems]{Beyond the PI Controllers in First-Order Time-Delay Systems}

\author{Gianpasquale Martelli}
\address{Via Domenico da Vespolate 8, 28079 Vespolate, Italy}
\curraddr{}
\email{gianpasqualemartelli@libero.it}
\thanks{}



\date{May 23, 2007}

\dedicatory{}

\begin{abstract}
In this paper the following three control systems for first-order time-delay plants are studied and compared: the feedback proportional-integral controller (PI), the Smith Predictor (SP) and a proposed variable structure consisting of two blocks. This structure acts as an open-loop proportional controller, after a setpoint change, and as a closed-loop integrating controller, when the error enters in a preset band. A chart, provided with the borderlines of the stability zone and with the curves of two design parameters, is implemented for each controller. The first parameter is the overshoot of the controlled variable, evaluated during a step change of the setpoint and made equal to a preset value. The second parameter, only for the PI and SP controllers, is the integral of the squared error (ISE), which  must have the minimum allowable value. The ISE is also assumed as  comparison index and the proposed controller appears as the best.
\end{abstract}

\maketitle

\section{Introduction}
First-order time-delay processes are examined in this paper, since  they are a simple but adequate approximation of most industrial plants. Moreover, controllers provided with two parameters are studied, since only systems, which can operate in a satisfying way when two performance indices are optimized, are taken into consideration.

The structure, represented in Fig. 1 and provided with a proportional-integral (PI) controller,  is the most used in time-delay systems thanks to its simplicity and to its robust performance in a wide range of operating conditions. Several tuning methods have been developed and their reliability and accuracy have been continuously increased. The first tuning rules, still prevalently implemented in the industrial systems  since very easy to use, have been proposed by Ziegler-Nichols \cite{bib1}. A better dynamical behaviour is assured by other methods, among them the performance criterion of Zhuan-Atherton \cite{bib2} and the deadbeat response rule \cite{bib3} deserve to be mentioned. Moreover, a new mathematical tool, consisting of the explicit solution of differential difference equations obtained by the method of steps \cite{bib4}, is recently available and therefore allows a more precise representation of the transient behaviour and an improved tuning of the control systems.

In the past some attempts of finding new control implementations, suitable  to eliminate some shortcomings occuring with the existing systems, have been made. The most significative is  the Smith Predictor (SP) \cite{bib5} (see Fig. 2), which, when  perfectly matched to the process,  has the very useful property to act as a PI controller, subject to a delayed setpoint but operating in a feedback system without delay.

In this paper a variable structure, depicted in Fig. 3 and based on already investigated and well-known theories and procedures, is presented. It indeed makes use of both open-loop proportional and closed-loop integrating control, of  the tracking function and of an adaptive technique updating the controller gain.
 
The selected tuning indices are the overshoot of the controlled variable, named $PO_{y}$, and the integral of the squared error ($ISE$), both evaluated during a  response to a step setpoint from the rated value to zero. It is also required an overshoot of the controller output, named $PO_{v}$, lower than a preset value, in order to avoid the controller windup. The $ISE$, calculated for a definite time  $t_{s}$ from the instant of the setpoint change, is then selected as comparison index. In order to reduce the computing complexity to an acceptable level, $t_{s}$ is assumed  equal to seven times the process time delay,  enough long to yield sound results. 

In Section 3 the new structure, described in Section 2, and the SP and PI controllers will be tuned and  compared. All the information needed for their implementation  will be included in charts, having as coordinates the controller parameters and provided with the stability zone borderlines and with the curves of the performance indices, whose analytical expressions are detailed in Appendices A, B and C. Finally, Section 4 contains some conclusive remarks. 

\begin{figure}[htbp]
\centering
\includegraphics{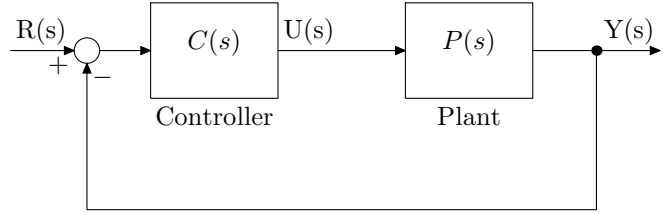}
\caption{PI controller}
\end{figure}

\begin{figure}[htbp]
\centering
\includegraphics{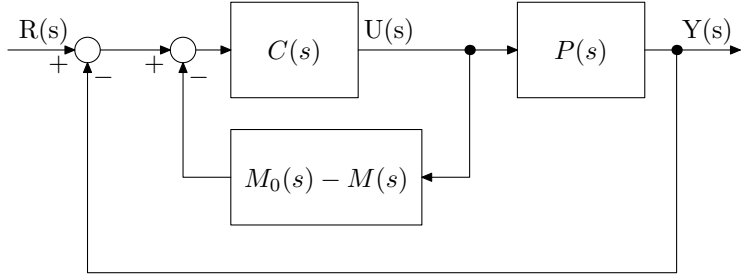}
\caption{SP controller}
\end{figure}

\section{Proposed controller}
The proposed controller, depicted in Fig. 3, consists of two blocks, $C_{1}$ and $C_{2}$, selected respectively for the first and the second operation mode. It is necessary to measure firstly the process gain, for example with an open-loop test, and to memorize it in the  controller.
The setpoint is continuously sampled and each value is compared with the previous ones, in order to check if its variation is inside a preset $ \pm B_{s}$ band (typically $B_{s}=2 \%$). When the first operation mode is working, the $C_{2}$ output tracks the $C_{1}$ output.

When a setpoint change bigger than $B_{s}$ occurs, the first operation mode, if not already working, is switched on, the controller acts as a open-loop proportional controller and the $C_{1}$ output becomes equal to the value corresponding to the actual setpoint.
When the error enters into the $ \pm B_{s}$ band, the second operation mode is switched on and the controller acts as a closed-loop integrator. 

While the second operation mode is on and until a new setpoint variation requiring the first operation mode occurs, an approach of the controlled variable to  the setpoint without oscillations outside the $ \pm B_{s}$ band is necessary. This can happen only if suitable values of  the $C_{2}$ parameter is chosen.
When the setpoint is constant or stays inside a narrow preset band for a preset time, the ratio between the controlled variable and the controller output is memorized in the block $C_{1}$ as updated process gain. 

\begin{figure}[htbp]
\centering
\includegraphics{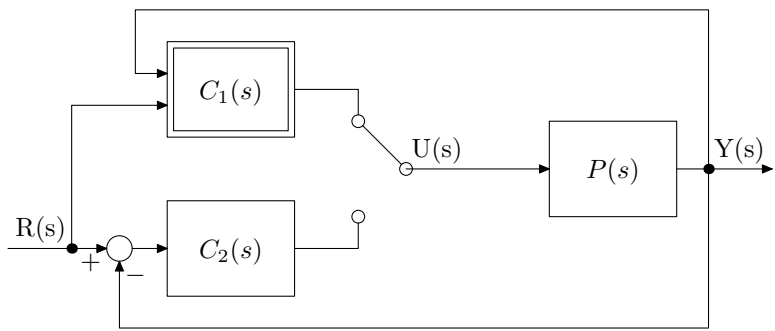}
\caption{Proposed controller}
\end{figure}

\section{Main results}
A very useful implementation tool is represented by the tuning chart, provided with the stability zone borderline and with the curves of the performance indices. The tuning point is simply the intersection of two curves, each referred to a preset value of a performance index; the stability of the system under examination is assured, if this point lies in the stability region. These charts are plotted in Figs. 4 (valid for $t_{p}=0.55$ as example), 5 and 6 and their coordinates are respectively $[h;h_{i}]$,$[h;h_{i}t_{p}]$ and $[t_{p};h_{i}]$. 

During the response to the setpoint change from the rated value to zero, the variables $r(t)$, $u(t)$ and $y(t)$, corresponding to $R(s)$, $U(s)$ and $Y(s)$, will be referred to their rated values.  The SP tuning indices $PO_{y}$ and $PO_{v}$  are directly evaluated by finding the minimum of analytical functions in the range from $t=0$ to $t= \infty$ and are given respectively by (\ref{eq:B.8}) and (\ref{eq:B.9}).  The index $PO_{y}$ of the other two controllers is the opposite of the minimum negative of 701 values of $y(t_{n})$ computed at equally spaced points; analogously  $PO_{v}$ is the minimum of  $v(t_{n})$. The used expressions are (\ref{eq:A.3}) and (\ref{eq:A.4}) for the PI and (\ref{eq:C.6}) or (\ref{eq:C.7}) and (\ref{eq:C.8}) for the proposed controller ($1<=n<=7$); the corresponding curves $\Gamma_{y}$ ($PO_{y}=0.0105$) and $\Gamma_{v}$ ($PO_{v}=0.1$) are plotted in the tuning charts (see Figs. 4 and 6).
The $ISE$ of each controller is evaluated with the triangle rule for the same time (from $t=0$ to $t=t_{s}=7$), according to  (\ref{eq:A.5}), (\ref{eq:B.10}) and (\ref{eq:C.9}).

\subsection{Transfer functions}
The transfer functions of the blocks included in Figs. 1, 2 and 3 are given by
\begin{equation}\label{eq:3.1}
P(s)= P_{0}(s) e^{-s \, L}= \frac{ K }{ 1+s \,T_{p}}\ e^{-s \, L}
\end{equation}
\begin{equation}\label{eq:3.2}
C(s)=K_{p}+ \frac{K_{i}}{s}\ 
\end{equation}
\begin{equation}\label{eq:3.3}
C_{2}(s)= \frac{K_{i}}{s}\ 
\end{equation}
where $K$ represents the plant steady-state gain,  $T_{p}$ the plant time constant, $L$ the plant time delay and $K_{p}$, $K_{i}$  the controller parameters.
It is convenient, in order to get equations independent of the real values of the parameters, to introduce the normalized time referred to the plant time delay $L$ ($ \sigma= s \,L$) and the dimensionless parameters $t_{p}=T_{p}/L$, $h=K \,K_{p}$ and $h_{i} = K \,K_{i}L$.  

Considering a matched SP ($M_{0}(\sigma)=P_{0}(\sigma)$ and $M(\sigma)=P(\sigma)$), the relationships between the controlled variable $Y(\sigma)$, the setpoint $R(\sigma)$ and the controller output $U(\sigma)$, based on the transfer functions from (\ref{eq:3.1}) to (\ref{eq:3.3}), are the following:
\begin{enumerate}

\item PI controller
\begin{equation}\label{eq:3.5}
Y(\sigma)=C(\sigma)P(\sigma)(R(\sigma)-Y(\sigma))=  \frac{h_{i}+  \sigma\,h}{\sigma(1+\sigma\,t_{p})}\ e^{-\sigma}(R(\sigma)-Y(\sigma))
\end{equation}
\begin{equation}\label{eq:3.6}
U(\sigma)=\frac{1}{P(\sigma)}\ Y(\sigma)= \frac{1+ \sigma \,t_{p}}{K}\ e^{ \sigma} Y(\sigma)
\end{equation}

\item SP controller
\begin{equation}\label{eq:3.7}
Y(\sigma)= \frac{C(\sigma)P_{0}(\sigma)}{1+C(\sigma)P_{0}(\sigma)}\ e^{- \sigma} R(\sigma)=  \frac{h_{i}+  \sigma\,h}{\sigma^{2}t_{p}+ \sigma (1+h)+h_{i}}\ e^{- \sigma}R(\sigma)
\end{equation}
\begin{equation}\label{eq:3.8}
U(\sigma)=\frac{1}{P(\sigma)}\ Y(\sigma)= \frac{1+ \sigma\,t_{p}}{K}\  Y(\sigma)
\end{equation}

\item proposed controller - first operation mode
\begin{equation}\label{eq:3.9}
Y(\sigma)=P(\sigma)U(\sigma)=  \frac{K e^{- \sigma}}{1+ \sigma \,t_{p} }\ U(\sigma)
\end{equation}

\item proposed controller - second operation mode
\begin{equation}\label{eq:3.10}
Y(\sigma)=C_{2}(\sigma)P(\sigma)(R(\sigma)-Y(\sigma))=  \frac{h_{i} }{\sigma (1+ \sigma \,t_{p}) }\ e^{- \sigma}(R(\sigma)-Y(\sigma))
\end{equation}
\begin{equation}\label{eq:3.11}
U(\sigma)=\frac{1}{P(\sigma)}\ Y(\sigma)= \frac{1+ \sigma \,t_{p}}{K}\ e^{ \sigma} Y(\sigma)
\end{equation}

\end{enumerate}

\subsection{Stability}
The stability of the SP controllers is theoretically assured for $1+h>0$, which always holds since $h>0$. Thr  instability may really happen if the process model transfer function is not correct, but in this paper this is excluded.
The following, proved in \cite{bib6} and confirmed in \cite{bib7}, hold for the PI ($h \neq 0$) and  the proposed ($h=0$) controllers:
\begin{enumerate}
\item The parameter $h$ must be included in the range from zero to $h_{u}$ given by
\begin{displaymath}
h_{u} = - \cos(z_{a})+t_{p}z_{a}sin(z_{a})
\end{displaymath}
where $z_{a}$ is the first positive solution of
\begin{displaymath}
tan(z_{a}) = - \frac{t_{p}}{1+t_{p}}\ z_{a}
\end{displaymath}
\item The parameters $h_{i}$, for a given $h$, must satisfy the following inequalities
 \begin{displaymath}
\delta_{r}(z_{1})<0 \, and \, \delta_{r}(z_{2})>0 \, and \, h_{i}>0
\end{displaymath}
 where 
\begin{displaymath}
\delta_{r}(z)= \frac{1}{L}\ [h_{i} - z \, sin(z) - t_{p}z^{2} cos(z)]
\end{displaymath}
and $z_{1}$ and $z_{2}$ are the first two positive roots of 
\begin{displaymath}
h+cos(z)-t_{p} z \, sin(z)=0
\end{displaymath}
\end{enumerate}

It is also useful to take into consideration  the stability phase margin (PM), whose expression is

\begin{displaymath}
\tan(z_{b}+PM)=-\frac{h_{i}+z_{b}^{2}h \,t_{p}}{z_{b}(h-h_{i}t_{p})}\
\end{displaymath}
where $z_{b}$ is the solution of

\begin{displaymath}
h^{2}+ \frac{h_{i}^{2}}{z_{b}^{2}}\ = 1+ t_{p}^{2}z_{b}^{2}.
\end{displaymath}
The curves corresponding to the stability region borderline and to $PM=30^{\circ}$, $PM=45^{\circ}$ and $PM=60 ^{\circ}$, named respectively  $\Gamma_{s}$, $\Gamma_{p1}$, $\Gamma_{p2}$ and $\Gamma_{p3}$, are plotted in the tuning charts (see Figs. 4 and 6). 

\subsection{Saturation}
One of the nonlinear issues of any controller is the actuator saturation, since from the practical point of view the actuator transfer function is bonded and the controller output may be out of the rated band.
This windup is avoided, in the simplest way, by switching off the controller integrating action when the actuator saturates, but this mode introduces some control discontinuities. A more sophisticated solution is given by the addition of the anti-windup scheme, depicted in Fig. 10.50 of  \cite{bib8}.
An other simple solution without discontinuities, adopted in this paper, consists in setting the rated value of the controller output lower than the extreme point of the linear part of the actuator transfer function and in tuning the controller in such a way that this point is never reached.

\subsection{Steadiness}
When the proposed controller is switched from the first to the second operation mode, the controlled variable, approaching the setpoint, should enter into and remain inside the preset $ \pm B_{s}$ band (typically $B_{s}=0.02$) until a new setpoint variation requiring the first operation mode occurs. This happens if the parameter $h_{i}$ is enough low and precisely if the tuning point is below the curve $\Gamma_{b}$, corresponding to $PO_{b}=B_{s}$, where $PO_{b}$ is the maximum absolute of 701 values of $y(t_{n})$, computed at equally spaced points with (\ref{eq:C.6}) or (\ref{eq:C.7}) (see fig. 6).

\subsection{Tuning}
\begin{enumerate}
\item PI controller\\
The tuning is performed according to the deadbeat response rule, presented in \cite{bib3} and summarized as follows: the  $ISE$ must be minimum and the controlled variable overshoot $PO_{y}$ must be equal to a preset value ($PO_{y}=0.0105$), provided that the controller output overshoot $PO_{v}$ is not higher than a fixed value ($PO_{v} \leq 0.10$).  
The tuning point $B$, having the minimum $ISE$ value and lying on $\Gamma_{y}$, is indicated in Fig. 4. 

\item SP controller\\
It is tuned with the same criterion as the PI controller.
Since the $PO_{v}$ and $ISE$ values on each point, lying on $\Gamma_{y}$, increases and decreases respectively when the coordinate $h$ increases, the intersection of $\Gamma_{y}$ and $\Gamma_{v}$, named $B$, must be assumed as tuning point. The damping borderline $\Gamma_{d}$, corresponding to $b=0$ according to (\ref{eq:B.4}), is also plotted in Fig. 5; the underdamped setpoint response happens for points above $\Gamma_{d}$.

\item proposed controller\\ 
This controller has one parameter ($h_{i}$) and, therefore, it is sufficient to require a controlled variable overshoot equal to a preset value ($PO_{y}=0.0105$).  The point, lying on the related curve $\Gamma_{y}$ in Fig.6 in correspondence of the coordinate $t_{p}$, can be considered as tuning point, since below both the controller output curve  $\Gamma_{v}$   ($PO_{v}=0.0105$) and the steadiness curve $\Gamma_{b}$ ($PO_{b}=B_{s}=0.02$). 
\end{enumerate}
The tuning parameters $h$ and $h_{i}$ and the performance indices $ISE$, related to some values of $t_{p}$ from $t_{p}=0.1$ to $t_{p}=10$, are gathered in Table \ref{Tab1} (as per Tables 1. and 2. of \cite{bib3} for the PI controller). It is plain that the proposed controller has the minimum values of the $ISE$ for all $t_{p}$ and therefore it can be elected the best.

\begin{figure}[htbp]
\centering
\includegraphics{martelliE1fig4.eps}
\caption{PI controller tuning chart for $t_{p}=0.55$} 
\end{figure}

\begin{figure}[htbp]
\centering
\includegraphics{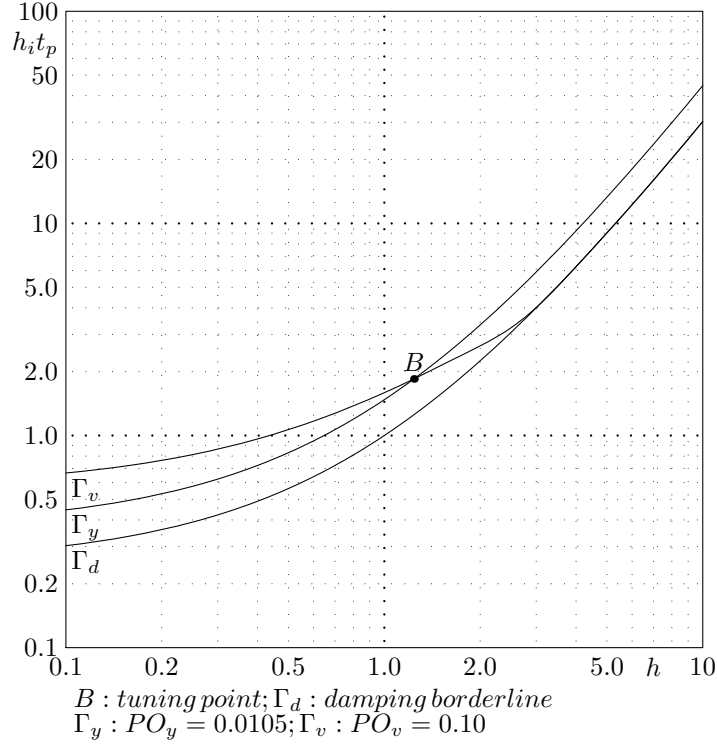}
\caption{Smith Predictor controller tuning chart}
\end{figure}

\begin{figure}[htbp]
\centering
\includegraphics{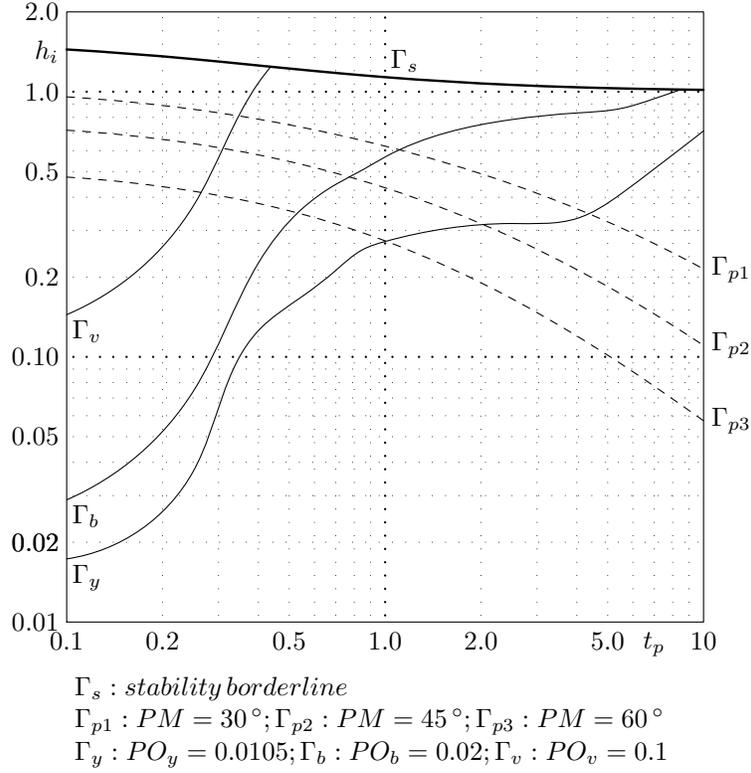}
\caption{Proposed controller tuning chart}
\end{figure}

\begin{table}[htbp]
\caption{Parameters and ISE of the controllers }
\label{Tab1}
\centering
\begin{tabular}{|r @{.} l|r @{.} l|r @{.} l|r @{.} l|r @{.} l|r @{.} l|r @{.} l|r @{.} l|r @{.} l|}
\hline
\multicolumn{2}{|c|}{} &
\multicolumn{2}{|c|}{PI} &
\multicolumn{2}{|c|}{PI} &
\multicolumn{2}{|c|}{SP} &
\multicolumn{2}{|c|}{SP} &
\multicolumn{2}{|c|}{Prop.} &
\multicolumn{2}{|c|}{PI} &
\multicolumn{2}{|c|}{SP} &
\multicolumn{2}{|c|}{Prop.} \\
\cline{3-18}
\multicolumn{2}{|c|}{$t_{p}$} &
\multicolumn{2}{|c|}{$h$} &
\multicolumn{2}{|c|}{$h_{i}$} &
\multicolumn{2}{|c|}{$h$} &
\multicolumn{2}{|c|}{$h_{i}$} &
\multicolumn{2}{|c|}{$h_{i}$} &
\multicolumn{2}{|c|}{ISE} &
\multicolumn{2}{|c|}{ISE} &
\multicolumn{2}{|c|}{ISE} \\
\hline
 0&10 & 0&45 & 0&787 & 1&239 & 18&490 & 0&017 
      & 1&524 & 1&083 & 1&051  \\
\hline
 0&25 & 0&50 & 0&738 & 1&239 & 7&396 & 0&037 
      & 1&674 & 1&207 & 1&125  \\
\hline 
 0&40 & 0&60 & 0&724 & 1&239 & 4&622 & 0&126 
      & 1&788 & 1&331 & 1&200  \\
\hline
 0&55 & 0&70 & 0&737 & 1&239 & 3&362 & 0&169 
      & 1&869 & 1&456 & 1&275  \\
\hline
 0&70 & 0&92 & 0&763 & 1&239 & 2&641 & 0&212 
      & 1&945 & 1&580 & 1&350  \\
\hline
 0&85 & 1&10 & 0&766 & 1&239 & 2&175 & 0&254 
      & 2&037 & 1&704 & 1&425  \\
\hline
 1&00 & 1&15 & 0&744 & 1&239 & 1&849 & 0&272 
      & 2&129 & 1&829 & 1&500  \\
\hline
 2&50 & 2&10 & 0&682 & 1&239 & 0&740 & 0&318 
      & 2&939 & 3&069 & 2&240  \\
\hline
 4&00 & 3&00 & 0&654 & 1&239 & 0&462 & 0&333 
      & 3&582 & 4&175 & 2&900  \\
\hline
 5&50 & 3&80 & 0&633 & 1&239 & 0&336 & 0&412 
      & 4&077 & 4&971 & 3&440  \\
\hline
 7&00 & 4&75 & 0&628 & 1&239 & 0&264 & 0&512 
      & 4&458 & 5&503 & 3&870  \\
\hline
 8&50 & 6&00 & 0&640 & 1&239 & 0&218 & 0&611 
      & 4&754 & 5&862 & 4&214  \\
\hline
 10&00 & 6&65 & 0&622 & 1&239 & 0&185 & 0&711 
      & 4&993 & 6&110 & 4&494  \\
\hline
\end{tabular}
\end{table}

\section{Conclusions}
In this paper a new controller is presented. Its main characteristics are:
\begin{enumerate}
\item during the switching from the first to the second operation mode, it is bumpless since the controller outputs of both blocks have the same value 
\item  it is adaptive, since the process gain is continuously monitored and accordingly updated, when a discrepancy appears between model and process values
\item  the allowable  values of the controller parameters, which avoid the actuator saturation, are easily deduced from the tuning chart
\item the process mismatch can be easily taken into consideration in the tuning chart, since the process parameters are included in their coordinates and suitable margins can be previously selected
\item it has zero steady-state error during the response to a step setpoint also for an  integrating process
\item  the lowest values of the $ISE$ can be achieved
\item it can be easily implemented in the modern digital controllers.
\end{enumerate}
The proposed controller may replace the PI controller, especially when the system is manually controlled, significative setpoint changes are needed and lowest $ISE$ values are required.

\appendix{}

\section{Response to a step setpoint of the PI controller}
During a response to a step setpoint from $r=1$ to $r=0$ applied to the closed-loop system in steady condition at $t=0$, the differential equation, related to the controlled variable $y$ and deduced from (\ref{eq:3.5}), is given by
\begin{equation}\label{eq:A.1}
\begin{split}
 \frac{dy(t)}{dt}\ &+ t_{p}  \frac {d^{2}y(t)}{dt^{2}}\ \\ &= (h_{i}+  h  \frac{d}{dt}\  ) (r(t-1)-y(t-1)) \quad  for \, t>=0
\end{split}
\end{equation}
subject to an initial condition of the form
\begin{equation}\label{eq:A.2}
 y(t)=1 \quad for \quad -1<=t<=0
\end{equation}
The analytical solution of (\ref{eq:A.1}), obtained with the method of steps in  \cite{bib4},  consists of a set of functions $y_{n}(t_{n})$, each valid for $n-1<t<n$ and given by 
\begin{equation}\label{eq:A.3}
\begin{split}
y_{1}(t_{1}) &= 1\\
y_{n}(t_{n}) &= \sum_{i=0}^{i=n-1}A_{n,i}t_{n}^{i}+ e^{-t_{n}/t_{p}} \sum_{j=0}^{j=n-2}B_{n,j}t_{n}^{j} \quad for \quad n>1
\end{split}
\end{equation}
where $t_{n}$ has the temporal origin at $t=n-1$.
In Example no.2 of  \cite{bib4} there are also the recursive expressions suitable for the evaluation of the coefficients $A_{n,i}$ and $B_{n,j}$.

Considering $t_{n}=t_{n+1}$, from (\ref{eq:3.6}) one obtains for $v_{n}(t_{n})=K \, u_{n}(t_{n})$ 
\begin{equation}\label{eq:A.4}
v_{n}(t_{n})=y_{n+1}(t_{n+1})+ t_{p} \frac{dy_{n+1}(t_{n+1})}{dt_{n+1}}\
\end{equation}
Finally, the integral of the squared error is given by
\begin{equation}\label{eq:A.5}
\begin{split}
ISE= &+ \frac{1}{2}\ 0.01 \,(y_{1}(0))^{2} - \frac{1}{2}\ 0.01 \, (y_{t_{s}}(1))^{2}\\ &+ 0.01 \sum_{n=1}^{n=t_{s}} \sum_{\tau=1}^{\tau=100} (y_{n} (0.01 \, \tau))^{2}\\ 
\end{split}
\end{equation}
where $t_{s}=7$.

\section{Response to a step setpoint of the SP controller}
During a response to a step setpoint from $r=1$ to $r=0$ applied to the closed-loop system in steady condition at $t=0$, the equation, related to the controlled variable $y$ and deduced from (\ref{eq:3.7}), and its underdamped solution are given by
\begin{equation}\label{eq:B.1}
\begin{split}
&t_{p} \frac{d^{2}y(t)} {dt^{2}}\ + (1+h) \frac{dy(t)} {dt}\ +h_{i} \, y(t)  = h_{i} \quad for \, 0<t<1\\ 
&t_{p} \frac{d^{2}y(t)} {dt^{2}}\ + (1+h) \frac{dy(t)} {dt}\ +h_{i} \, y(t)  = 0 \quad for \, t>1 
\end{split} 
\end{equation}
\begin{equation}\label{eq:B.2}
\begin{split}
y(t)&=1 \quad for \quad 0<t<1 \\
y(t)&=e^{-a \,(t-1)}(cos(b \, (t-1))+ \frac {a}{b}\ sin(b \,(t-1))) \quad for \quad t>1 
\end{split}
\end{equation}
where
\begin{displaymath}
a=\frac{0.5(1+h) }{t_{p}}\
\end{displaymath}
\begin{equation}\label{eq:B.4}
b= \frac{0.5}{t_{p}}\ \left(4 \, h_{i} \, t_{p}-(1+h)^{2} \right)^{0.5}
\end{equation}
From (\ref{eq:3.8}) and (\ref{eq:B.2}) one obtains for $v(t)=K \, u(t)$
\begin{equation}\label{eq:B.5}
\begin{split}
v(t)&=y(t)+ t_{p}\frac {dy(t)}{dt}\ \\ &=e^{-a \,(t-1)}(cos(b \, (t-1))- \frac {-a+t_{p}(a^{2} +b^{2})}{b}\ sin(b \,( t-1)))
\end{split}
\end{equation}
The derivatives of $y(t)$ and $v(t)$ with respect to $t$ are given respectively by
\begin{equation}\label{eq:B.6}
\frac{dy(t)}{dt}\ = - e^{-a \,(t-1)} \frac {a^{2}+b^{2}}{b}\ sin(b \,(t-1))
\end{equation}
\begin{equation}\label{eq:B.7}
\frac{dv(t)}{dt}\ = e^{-a \,(t-1)}(a^{2}+b^{2})(-t_{p}cos(b \, (t-1))+ \frac {a \,t_{p}-1}{b}\ sin(b \, (t-1)))
\end{equation}
Introducing the solution $b(t-1)$ of (\ref{eq:B.6}), equated to zero, in (\ref{eq:B.2}), we obtain for the overshoot $PO_{y}$ 
\begin{equation}\label{eq:B.8}
PO_{y} = -e^{- \pi \,a/b}
\end{equation}
Analogously introducing the solution $b(t-1)$ of (\ref{eq:B.7}), equated to zero, in (\ref{eq:B.5}), we obtain for the overshoot $PO_{v}$ 
\begin{equation}\label{eq:B.9}
PO_{v} =- e^{ - \phi \, a/b}(1-2 \,a \,t_{p}+t_{p}^{2}(a^{2}+b^{2}))^{0.5}
\end{equation}
where $\phi$ is the first positive root of $tan(\phi)=b \,t_{p}/(a \,t_{p}-1)$.
Finally, the integral of the squared error is given by
\begin{equation}\label{eq:B.10}
ISE =1+ ISE_{a}(t_{s})-ISE_{a}(1)
\end{equation}
where $t_{s}=7$ and
\begin{displaymath}
\begin{split}
ISE_{a}(t) = &\frac{e^{-2 \,a(t-1)}}{4 \, a \,b^{2}(a^{2}+b^{2})}\  ( -(a^{2}+b^{2})^{2}\\ &+a^{2}(a^{2}-3 \,b^{2})cos(2 \,b(t-1))+a \,b (-3 \,a^{2}+b^{2})sin(2 \,b(t-1)) )
\end{split}
\end{displaymath}

\section{Response to a step setpoint of the proposed controller}
During a response to a step setpoint from $r=1$ to $r=0$ applied to the closed-loop system in steady condition at $t=0$, the following holds:

\begin{enumerate}
\item First operation mode\\
From  (\ref{eq:3.9}) one obtains
\begin{equation}\label{eq:C.1}
\begin{split}
&y_{a}(t)=1 \quad for \, 0<=t<=1\\
&y_{b}(t)=e^{(-t+1)/t_{p}} \quad  for \, t>=1
\end{split}
\end{equation}
since $y_{b}(t)$ is the solution of
\begin{displaymath}
y_{b}(t)+t_{p}\frac{dy(t)}{dt}\ =0
\end{displaymath}
subject to the initial condition $y_{b}(1)=y_{a}(1)=1$.
Considering the identity $y_{b}(1+t_{q})=B_{s}$, the time $t_{q}$ required by $y_{b}$ to reach $B_{s}$ is given by $t_{q}=t_{p}ln(1/B_{s})$. 
\item Second operation mode\\ 
From  (\ref{eq:3.10}) one obtains
\begin{equation}\label{eq:C.3}
\begin{split}
 \frac{dy(t)}{dt}\ &+ t_{p}  \frac {d^{2}y(t)}{dt^{2}}\ \\ &= (h_{i}+  h  \frac{d}{dt}\  ) (r(t-1)-y(t-1)) \quad  for \, t>=1+t_{q}
\end{split}
\end{equation}
subject to the initial condition
\begin{enumerate}
\item $t_{q}<1$
\begin{displaymath}
\begin{split}
y(t) &=y_{a}(t)=1 \quad for \quad t_{q}<t<1\\
y(t) &=y_{b}(t) =e^{(-t+1)/t_{p}} \quad for \quad 1<t<1+t_{q}
\end{split}
\end{displaymath}
\item $t_{q}>1$
\begin{displaymath}
y(t)= y_{b}(t)=e^{(-t+1)/t_{p}} \quad for \quad t_{q}<t<1+t_{q}
\end{displaymath}
\end{enumerate}

The analytical solution of (\ref{eq:C.3}), obtained with the method of steps in  \cite{bib4},  consists of a set of functions $y_{n}(t_{n})$, each valid for $n>=1$ and given by 

\begin{enumerate}
\item $t_{q}<1$
\begin{equation}\label{eq:C.6}
\begin{split}
y_{n,1}(t_{n}) &= \sum_{i=0}^{i=n}A_{n,1,i}t_{n}^{i}+ e^{-t_{n}/t_{p}} \sum_{j=0}^{j=n-1}B_{n,1,j}t_{n}^{j} \quad for \quad 0<t_{n}<1-t_{q}\\
y_{n,2}(t_{n}) &= \sum_{i=0}^{i=n-1}A_{n,2,i}t_{n}^{i}+ e^{-t_{n}/t_{p}} \sum_{j=0}^{j=n}B_{n,2,j}t_{n}^{j} \quad for \quad 1-t_{q}<t_{n}<1
\end{split}
\end{equation}
\item $t_{q}>1$
\begin{equation}\label{eq:C.7}
y_{n}(t_{n}) = \sum_{i=0}^{i=n-1}A_{n,i}t_{n}^{i}+ e^{-t_{n}/t_{p}} \sum_{j=0}^{j=n}B_{n,j}t_{n}^{j} \quad for \quad 0<t_{n}<1
\end{equation}
\end{enumerate}
where $t_{n}$ has the temporal origin at $t=t_{q}+n$.
The initial conditions are given by
\begin{enumerate}
\item $t_{q}<1$: $y_{0,1}(t_{0})=1$ and $y_{0,2}(t_{0}) = e^{-t_{0}/t_{p}} \, e^{(-t_{q}+1)/t_{p}}$
\item $t_{q}>1$: $y_{0}(t_{0})=e^{-t_{0}/t_{p}} \,e^{(-t_{q}+1)/t_{p}}$
\end{enumerate}
In Example no.1 of  \cite{bib4} there are also the recursive expressions suitable for the evaluation of the coefficients $A_{n,i}$, $B_{n,j}$, $A_{n,1,i}$, $B_{n,1,j}$, $A_{n,2,i}$ and $B_{n,2,j}$.
\end{enumerate}
Considering $t_{n}=t_{n+1}$, from (\ref{eq:3.11}) one obtains for $v_{n}(t_{n})=K \, u_{n}(t_{n})$ 
\begin{equation}\label{eq:C.8}
v_{n}(t_{n})=y_{n+1}(t_{n+1})+ t_{p} \frac{dy_{n+1}(t_{n+1})}{dt_{n+1}}\
\end{equation}
Finally, the integral of the squared error is given by
\begin{equation}\label{eq:C.9}
\begin{split}
ISE= &1 + \frac{t_{p}}{2}\ (1-e^{-2 \,t_{q}/t_{p}}) +ISE_{b} \quad for \, 0<t_{q}<t_{s}-1 \\
ISE= &1 + \frac{t_{p}}{2}\ (1-e^{-2(t_{s} -1)/t_{p}})
\quad for \, t_{s}-1<t_{q} 
\end{split}
\end{equation}
where $t_{s}=7$ and $ISE_{b}$, evaluated from $t=1+t_{q}$ to $t=t_{s}$ with an expression similar to (\ref{eq:A.5}), is the contribute of the second operation mode.

\end{document}